\documentclass[dvipsnames]{daj} 

\usepackage{amsthm,amssymb,amsmath,tikz}
\usetikzlibrary{patterns}
\numberwithin{equation}{section}
\newtheorem{theorem}{Theorem}[section]
\newtheorem{lemma}[theorem]{Lemma}
\theoremstyle{remark}
\newtheorem{remark}[theorem]{Remark}
\theoremstyle{definition}
\newtheorem{definition}[theorem]{Definition}
\newtheorem{example}[theorem]{Example}

\dajAUTHORdetails{%
  title = {Oscillation Estimates of Eigenfunctions via the Combinatorics of Noncrossing Partitions}, 
  author = {Vera Mikyoung Hur, Mathew A. Johnson, and Jeremy L. Martin},
  plaintextauthor = {Vera Mikyoung Hur, Mathew A. Johnson, Jeremy L. Martin},
  keywords = {fractional Schr\"{o}dinger; eigenfunction; oscillation; noncrossing partition},
}   

\dajEDITORdetails{%
   year={2017},
   number={13},
   received={9 September 2016},   
   revised={10 February 2017},    
   published={5 September 2017},  
   doi={10.19086/da.2102},       
}   

\begin{document}

\begin{frontmatter}[classification=text]

\title{Oscillation Estimates of Eigenfunctions via the Combinatorics of Noncrossing Partitions} 

\author[vmh]{Vera Mikyoung Hur\thanks{Supported by the National Science Foundation under the grant CAREER DMS-1352597, an Alfred P. Sloan Research Fellowship, a Simons Fellowship in Mathematics, and the University of Illinois at Urbana-Champaign under the Arnold O. Beckman Research Award RB14100 and a Beckman Fellowship.}}
\author[maj]{Mathew A. Johnson\thanks{Supported by the National Science Foundation under grants DMS-1614785 and DMS-1211183.}}
\author[jlm]{Jeremy L. Martin\thanks{Supported by Simons Foundation Collaboration Grant number 315347.}}

\begin{abstract}
We study oscillations in the eigenfunctions for a fractional Schr\"odinger operator on the real line. 
An argument in the spirit of Courant's nodal domain theorem applies to an associated local problem in the upper half plane and provides a bound on the number of nodal domains for the extensions of the eigenfunctions. 
Using the combinatorial properties of noncrossing partitions, we turn the nodal domain bound into an estimate for the number of sign changes in the eigenfunctions. We discuss applications in the periodic setting and the Steklov problem on planar domains.
\end{abstract}
\end{frontmatter}



\section{Introduction}\label{s:intro}

The purpose of this note is to estimate the number of sign changes in the eigenfunctions for a \emph{fractional 
Schr\"odinger operator}
\[
H=\left(-\frac{d^2}{dx^2}\right)^{\alpha/2}+V(x),
\]
where for $0<\alpha<2$, the fractional Laplacian $(-d^2/dx^2)^{\alpha/2}$ is defined via its multiplier $|\xi|^\alpha$ in the Fourier space, and $V$ is an appropriate potential. Throughout, we assume that $V$ is real valued and bounded. 

Fractional Laplacians arise in a variety of applications, including quantum mechanics, minimal surfaces, phase transition, crystal dislocation, and anomalous diffusion. 
We are motivated by problems related to wave motion. For instance, consider the fractional Korteweg-de Vries (fKdV) equation
\[
u_t+u_x+\left(-\frac{\partial^2}{\partial x^2}\right)^{\alpha/2}u_x+f(u)_x=0,
\]
the fractional Benjamin--Bona--Mahony (fBBM) equation
\[
u_t+\left(-\frac{\partial^2}{\partial x^2}\right)^{\alpha/2}u_t+u_x+f(u)_x=0,
\]
and the fractional nonlinear Schr\"odinger (fNLS) equation
\[
iu_t-\left(-\frac{\partial^2}{\partial x^2}\right)^{\alpha/2}u+f(|u|)u=0
\]
for $0<\alpha<2$ and an appropriate $f$. Here the fractional Laplacian comes from modeling dispersion. For (fKdV) and (fBBM), we assume that $u=u(x,t)$ is real valued; for (fNLS), $u$ may be complex valued.  
In many examples of interest, the existence of solitary and periodic traveling waves of (fKdV), (fBBM), and (fNLS)  follows from variational or perturbative arguments; see, e.g., \cite{BH14, HJ15, LPS15, A16, KSM14, HP16, CJ16} and references therein.  Linearizing  (fKdV), (fBBM), or (fNLS) about a traveling wave, one arrives at 
an operator of the form $(-d^2/dx^2)^{\alpha/2}+V(x)$, where $V$ depends on the underlying traveling wave.
The spectral properties of such a fractional Schr\"odinger operator, including the number of negative eigenvalues, play a central role in the study of the stability and instability of the underlying traveling wave, and for several other purposes; see, e.g., \cite{Lin08, LPS15, BH14, HJ15, J13, KS14, HP16, CJ16}.

When $\alpha=2$, so that $H=-d^2/dx^2+V(x)$ is the classical Schr\"odinger operator, one may use ordinary differential equations (ODE) techniques to locate the eigenvalues and count the number of zeros in the eigenfunctions.

\begin{theorem}[Sturm's oscillation theorem]\label{T:sturm}
Let $\alpha=2$, $V:\mathbb{R}\to\mathbb{R}$ and $V\in L^1(\mathbb{R})\bigcap L^\infty(\mathbb{R})$. 
For an integer $N\geq 1$, suppose that the $L^2(\mathbb{R})$ spectrum of $H=-d^2/dx^2+V(x)$ contains
at least $N$ eigenvalues 
\begin{equation}\label{E:eigen2}
\lambda_1<\lambda_2<\cdots<\lambda_N<0=\min\sigma_{\rm ess}(H).
\end{equation}
Then for each integer $n=1,2,\ldots N$, a real and $C^0(\mathbb{R})$ eigenfunction associated with the 
eigenvalue $\lambda_n$ has exactly $n-1$ zeros in $\mathbb{R}$.
\end{theorem}

See, e.g., \cite{KPbook} for details.  We emphasize that the strict inequalities in \eqref{E:eigen2} imply that each eigenvalue is simple and the associated eigenspace is one dimensional. 

Theorem~\ref{T:sturm} counts the exact number of zeros in the eigenfunctions in one dimension.  
Courant's nodal domain theorem (see, e.g., \cite{CH}) extends this to higher dimensions and offers an upper bound on the number of nodal domains (see \eqref{def:nodal domain}) of the eigenfunctions. 
However, the eigenvalues need not be simple and, furthermore, the nodal domain bound is not in general sharp.

When $0<\alpha<2$, so that $H=(-d^2/dx^2)^{\alpha/2}+V(x)$ is a \emph{nonlocal} operator, 
ODE techniques may not directly apply.
Rather, one needs to develop appropriate substitutes based on different methods.

Recently, Frank and Lenzmann \cite{FL13} successfully obtained a sharp bound on the number of sign changes in the second eigenfunction of $H=(-d^2/dx^2)^{\alpha/2}+V(x)$ for $0<\alpha<2$. 
Specifically, they recalled that the fractional Laplacian may be regarded as the Dirichlet-to-Neumann operator for an appropriate local problem in the upper half plane (see, e.g., \cite{CS07, GZ03}), whereby they developed a variational principle for the eigenvalues and eigenfunctions of $H$. An argument in the spirit of Courant's nodal domain theorem then applies and provides an upper bound on the number of nodal domains in the upper half plane for the extensions of the eigenfunctions. They made a topological argument to turn the nodal domain bound into the number of sign changes in the second eigenfunction; see Section~\ref{s:analysis} for some details. 
Recently, two of the authors \cite{HJ15,CJ16} extended the argument to periodic potentials, subject to either periodic or anti-periodic boundary conditions. 

Unfortunately, it is not clear how to extend the topological argument in \cite{FL13} and \cite{HJ15,CJ16} to higher eigenfunctions without recourse to a tedious case-by-case study. A satisfactory oscillation theory of eigenfunctions therefore seems lacking for the fractional Schr\"odinger operators and other nonlocal operators. Here we resort to the combinatorial properties of \emph{noncrossing partitions} to estimate the number of sign changes in the higher eigenfunctions of $H=(-d^2/dx^2)^{\alpha/2}+V(x)$ for $0<\alpha<2$. Thereby, we generalize Theorem~\ref{T:sturm} to fractional Schr\"odinger operators and the result in \cite{FL13} to all eigenfunctions.

Introduced by Kreweras \cite{Kreweras}, noncrossing partitions are a fundamental tool in algebraic combinatorics, with applications in numerous branches of mathematics, ranging from representation theory to free probability, and in other sciences, such as molecular biology. We encourage the interested reader to \cite{Simion, Armstrong} and references therein. Here we add a new application of noncrossing partitions to spectral theory. Briefly speaking, a noncrossing partition decomposes a finite and totally ordered set into pairwise disjoint subsets, called \emph{blocks}, which may be geometrically represented by polygons with vertices on a circle; see, e.g., Figure~\ref{F:ncp}. In the present setting, the vertices correspond to the points in $\mathbb{R}$ where the sign of the eigenfunction is either positive or negative. The points in each block are contained in the same nodal domain, and the noncrossing condition means that the nodal domains are disjoint. We devise a combinatorial argument to establish a relationship between the number of nodal domains of a continuous function in a planer domain and the number of sign changes of the function on the boundary.

We emphasize that since we use an argument in the spirit of Courant's nodal domain
theorem, the result is not sharp and should not be expected to be.
On the other hand, our approach is based on a general combinatorial argument that ought to be readily adapted to other related problems. To illustrate this, we apply our argument to the periodic setting as well as to the Steklov problem on bounded domains in $\mathbb{R}^2$ of arbitrary genus.  Specific features of each problem may help to improve the result, but we do not pursue the direction here.  Nevertheless, we demonstrate by example that our result is sharp for the Steklov problem. 

The organization of the article is as follows. In Section~\ref{s:analysis}, we recall some known facts about $H=(-d^2/dx^2)^{\alpha/2}+V(x)$, acting on $L^2(\mathbb{R})$, for $0<\alpha<2$.
In Section~\ref{s:comb}, we develop the combinatorics of noncrossing partitions, assuming no prior familiarity with the subject on the reader's part. In Section~\ref{s:application}, we combine the results in Section~\ref{s:analysis} and Section~\ref{s:comb} to estimate the number of sign changes in the eigenfunctions of $H$, and we discuss applications to the periodic setting as well as the Steklov problem.

\section{Nodal domain count}\label{s:analysis}

Let $0<\alpha<2$, $V:\mathbb{R} \to\mathbb{R}$ and $V\in L^\infty(\mathbb{R})$. (One may relax this to $L^p(\mathbb{R})$ for some $p\geq 1$; see, e.g., \cite{CMS} for details.) Note that 
\[
H=\left(-\frac{d^2}{dx^2}\right)^{\alpha/2}+V(x)
\] 
is a self-adjoint operator on $L^2(\mathbb{R})$ with domain $H^{\alpha}(\mathbb{R})$, and its essential
spectrum is $\sigma_{\rm ess}(H)=[b,\infty)$ for some $b\in\mathbb{R}$. (We may assume $b=0$ after adding an appropriate constant to $V$.)  Since $H$ is a real operator, we may assume that an eigenfunction of $H$ is real valued. Moreover, note that any $L^2(\mathbb{R})$ eigenfunction of $H$ is continuous 
and bounded over $\mathbb{R}$; see, e.g., \cite{CMS}.

Here we recall some known facts about $H=(-d^2/dx^2)^{\alpha/2}+V(x)$. Details may be found in \cite[Section~3]{FL13}. We merely touch on the main points.

We begin by observing that the fractional Laplacian $(-d^2/dx^2)^{\alpha/2}$ may be interpreted as the Dirichlet-to-Neumann operator for an appropriate local problem in the upper half plane $\mathbb{P}^2:=\{(x,y)\in\mathbb{R}^2:y>0\}$; see, e.g., \cite{CS07, GZ03} for details. Specifically, consider the boundary value problem
\begin{equation}\label{E:extension}
\begin{cases}
\nabla\cdot(y^{1-\alpha}\nabla w)=0 \qquad &\textrm{for}\quad(x,y)\in\mathbb{P}^2,\\
w=f&\textrm{on}\quad\partial\mathbb{P}^2=\mathbb{R}\times\{0\},
\end{cases}
\end{equation}
whose solution may be written explicitly as
\[
w(x,y)=\int_{\mathbb{R}}\frac{1}{y}P\left(\frac{x-z}{y}\right)f(z)~dz,\qquad\text{where}\quad 
P(x)=\frac{1}{\left(1+x^2\right)^{(1+\alpha)/2}}
\]
up to multiplication by a constant. Let 
\begin{equation}\label{def:E}
E_\alpha(f)=w,
\end{equation}
and we shall drop the subscript $\alpha$ when there is no ambiguity. For any $f\in H^{\alpha}(\mathbb{R})$, note that
\begin{equation}\label{e:limit}
\lim_{y\to 0+} y^{1-\alpha}\frac{\partial E_\alpha(f)}{\partial y}(\cdot,y)=-C_\alpha\left(-\frac{d^2}{dx^2}\right)^{\alpha/2} f
\qquad\text{and}\qquad\lim_{y\to 0+}E_\alpha(f)(\cdot,y)=f
\end{equation}
in appropriate functional spaces for some explicit constant $C_\alpha>0$. 

To proceed, one may characterize the eigenvalues and eigenfunctions of $H=(-d^2/dx^2)^{\alpha/2}+V(x)$ in terms of the Dirichlet type functional 
\begin{equation}\label{def:energy}
\iint_{\mathbb{P}^2}y^{1-\alpha}|\nabla w(x,y)|^2~dx\:dy+\int_{\mathbb{R}}V(x)|w(x,0)|^2~dx,
\end{equation}
which is defined in an appropriate function space for $w=w(x,y)$ in $\mathbb{P}^2$, where $w(x,0)$ denotes the trace of $w(x,y)$ on the boundary $\partial\mathbb{P}^2=\mathbb{R}\times \{0\}$.  Specifically, for an integer $N\geq 1$, suppose that the $L^2(\mathbb{R})$ spectrum of $H$ contains at least $N$ eigenvalues
\[
\lambda_1\leq \cdots \leq \lambda_N<\min \sigma_{\text ess}(H).
\]
Then for each $n=1,2,\dots, N$, $\lambda_n$ arises as the infimum of \eqref{def:energy} over an appropriate function space, subject to the condition that the trace on the boundary $\partial\mathbb{P}^2=\mathbb{R}\times\{0\}$ is orthogonal (with respect to the $L^2(\mathbb{R})$ inner product) to the $(n-1)$-dimensional subspace of $L^2(\mathbb{R})$ spanned by the eigenfunctions of $H$ associated with the eigenvalues $\lambda_k$ for $k=1,2,\dots, n-1$. 
Moreover, the infimum is achieved if and only if $w=E(f)$, where $f$ is a linear combination of the eigenfunctions of $H$ associated with $\lambda_n$. Indeed, for any $f\in H^\alpha(\mathbb{R})$ and $w=E(f)$, it follows that
\begin{align*}
C_\alpha^{-1}\iint_{\mathbb{P}^2}y^{1-\alpha}|\nabla w(x,y)|^2~dx\:dy
&=C_\alpha^{-1}\lim_{\epsilon\to 0^+}\iint_{\mathbb{R}\times(\epsilon,\infty)}y^{1-\alpha}|\nabla w(x,y)|^2~dx\:dy\\
&=-C_\alpha^{-1}\lim_{\epsilon\to 0^+}\epsilon^{1-\alpha}\int_{\mathbb{R}}\overline{\frac{\partial w}{\partial y}(x,y)}w(x,y)~dx\\
&=\int_{\mathbb{R}}\overline{(-d^2/dx^2)^{\alpha/2} f(x)}f(x)~dx \\
&=\left\|\left(-\frac{d^2}{dx^2}\right)^{\alpha/4} f\right\|_{L^2(\mathbb{R})}^2.
\end{align*}
Note that the second equality uses the divergence theorem and the third equality uses \eqref{e:limit}.

To continue, recall that any $L^2(\mathbb{R})$ eigenfunction $\phi$ of $H=(-d^2/dx^2)^{\alpha/2}+V(x)$ is real valued, bounded and continuous. Hence its extension $E(\phi)$ to the upper half plane is real valued and belongs to $C^0(\overline{\mathbb{P}^2})$. By the \emph{nodal domains} of $E(\phi)$, we mean the connected components in $\mathbb{P}^2$ of
\begin{equation}\label{def:nodal domain}
\mathbb{P}^2-\overline{\{(x,y)\in\mathbb{P}^2:E(\phi)(x,y)=0\}}.
\end{equation}
An argument in the spirit of Courant's nodal domain theorem offers an estimate on the number of nodal domains. 

\begin{theorem}[Courant's nodal domain theorem]\label{L:courant}
Let $0<\alpha<2$, $V:\mathbb{R}\to\mathbb{R}$ and $V\in L^\infty(\mathbb{R})$. For an integer $N\geq 1$, 
suppose that the $L^2(\mathbb{R})$ spectrum of $H=(-d^2/dx^2)^{\alpha/2}+V(x)$ 
contains at least $N$ eigenvalues
\begin{equation}
\lambda_1<\lambda_2\leq\ldots\leq\lambda_N<\min\sigma_{\rm ess}(H).
\end{equation}\label{E:eigen-alpha}
For each integer $n=1,2,\dots, N$, if $\phi_n\in H^{\alpha/2}(\mathbb{R})\bigcap C^0(\mathbb{R})$ is a real eigenfunction of $H$ associated with the eigenvalue $\lambda_n$ then its extension $E(\phi_n)$ has at most $n$ nodal domains in $\mathbb{P}^2$.
\end{theorem}

The proof uses the variational principle for $H$ (see \eqref{def:energy}) and the unique continuation of solutions of \eqref{E:extension}. See \cite{FL13} for details.

\begin{remark}\label{R:mult}
A remark is in order concerning how multiplicities are handled in Theorem~\ref{L:courant} and throughout.
The lowest eigenvalue $\lambda_1$ of $H=(-d^2/dx^2)^{\alpha/2}+V(x)$ is necessarily simple (see, e.g., \cite{FL13}), 
which we indicate by the strict inequality in \eqref{E:eigen-alpha}. 
But higher eigenvalues of $H$ are not guaranteed to be simple.
To illustrate the effects of multiplicities, suppose, for instance, that $\lambda_2$ is an eigenvalue
with algebraic multiplicity two, so that
\[
\lambda_1<\lambda_2=\lambda_3<\lambda_4\leq\ldots\leq\lambda_N.
\]
Theorem~\ref{L:courant} implies that if $\phi_2$ is an eigenfunction associated with $\lambda_2=\lambda_3$
then $E(\phi_2)$ has at most two nodal domains in $\mathbb{P}^2$, while if $\phi_4$ is an eigenfunction
associated with $\lambda_4$ then $E(\phi_4)$ has at most four nodal domains in $\mathbb{P}^2$.
Throughout, the oscillation results in the presence of multiplicities
must be interpreted similarly.
\end{remark}

\begin{definition}
For an integer $N\geq1$, we say that $\phi\in C^0(\mathbb{R})$ \emph{changes its sign $N$ times} in $\mathbb{R}$  if there exist $N+1$ points \[
x_1<x_2<\ldots<x_{N+1}
\]
such that $\phi(x_k)\neq 0$ for $k=1,2,\dots, N$ and $\text{sgn}\,(\phi(x_k))=-\text{sgn}\,(\phi(x_{k+1}))$ for $k=1,2,\dots, N$.
\end{definition}

We wish to turn the nodal domain bound in Theorem~\ref{L:courant} into a bound on the number of sign changes in the eigenfunctions of $H=(-d^2/dx^2)^{\alpha/2}+V(x)$. In other words, we want to relate the number of nodal domains of a function in $\mathbb{P}^2$ with the number of sign changes on the boundary $\partial\mathbb{P}^2=\mathbb{R}\times\{0\}$. 

Clearly, an eigenfunction associated with the lowest eigenvalue $\lambda_1$ does not change its sign. Otherwise, by continuity, its extension in the upper half plane would have at least two nodal domains, contrary to Theorem~\ref{L:courant}; see also Remark~\ref{R:sharp1} below. 
Frank and Lenzmann \cite{FL13} took matters further and established a sharp bound on the number of sign changes in an eigenfunction associated with the second eigenvalue $\lambda_2$.

\begin{theorem}[Oscillation in a second eigenfunction]\label{L:2nd}
Under the hypotheses of Theorem~\ref{L:courant}, an eigenfunction associated with $\lambda_2$ changes its sign at most twice in $\mathbb{R}$.
\end{theorem}

\begin{figure}[h]
\begin{center}
\begin{tikzpicture}
\newcommand{\crad}{0.7ex} 
\newcommand{\lht}{-.4} 
\newcommand{\xone}{-1.2}
\newcommand{\yone}{0.6}
\newcommand{\xtwo}{1.8}
\newcommand{\ytwo}{4.4}
\draw [very thick, <->] (-2.5,0) -- (5.5,0);
\draw [very thick, ->] (0,0) -- (0,2.5);
\draw[very thick, ProcessBlue] (\xtwo,0) arc (0:180:\xtwo/2-\xone/2);
\draw[very thick, Red] (\ytwo,0) arc (0:180:\ytwo/2-\yone/2);
\draw[thick, fill=Cyan] (\xone,0) circle (\crad);
\draw[thick, fill=Red] (\yone,0) circle (\crad);
\draw[thick, fill=Cyan] (\xtwo,0) circle (\crad);
\draw[thick, fill=Red] (\ytwo,0) circle (\crad);
\node at (\xone,\lht) {$x_1^+$};
\node at (\yone,\lht) {$x_1^-$};
\node at (\xtwo,\lht) {$x_2^+$};
\node at (\ytwo,\lht) {$x_2^-$};
\node at (-.7,1.4) {$\gamma^+$};
\node at (4.4,1.4) {$\gamma^-$};
\end{tikzpicture}
\end{center}
\caption{Graphical depiction of the proof of Theorem~\ref{L:2nd}.\label{F:twocrossing}}
\end{figure}
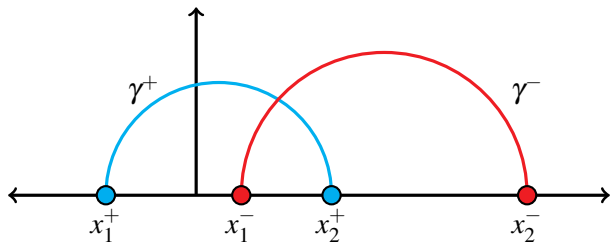

\begin{proof}
The proof is in \cite{FL13}. Here we sketch some details to aid the subsequent discussion.

Let $\phi_2$ denote an eigenfunction associated with $\lambda_2$. Suppose on the contrary that $\phi_2$ changes its sign at least three times in $\mathbb{R}$. Then there must exist four real numbers
\[
x^+_1<x^-_1<x^+_2<x^-_2
\]
such that, without loss of generality,
\[
\phi_2(x^+_k)>0\quad\text{and}\quad \phi_2(x^-_k)<0 \qquad\text{for $k=1,2$.}
\] 
Theorem~\ref{L:courant} asserts that the extension $E(\phi_2)$ to the upper half plane has exactly two nodal domains, one where $E(\phi_2)$ is positive and the other where it is negative. On the other hand, we may find a continuous curve $\gamma^+$ in $\mathbb{P}^2$ connecting $x^+_1$ and $x^+_2$ in the positive nodal domain and, likewise, a continuous curve $\gamma^-$ connecting $x^-_1$ and $x^-_2$ in the negative nodal domain; see Figure~\ref{F:twocrossing}. The Jordan curve theorem then dictates that $\gamma^+$ and $\gamma^-$ intersect in $\mathbb{P}^2$, which is a contradiction since nodal domains are disjoint by definition.
\end{proof}

\begin{remark}\label{R:sharp1}
Note that the lowest eigenvalue $\lambda_1$ of $H=(-d^2/dx^2)^{\alpha/2}+V(x)$ is simple 
and an associated eigenfunction is sign definite over $\mathbb{R}$. 
Since $H$ is self adjoint, it is easy to verify that 
each eigenfunction associated with an eigenvalue $\lambda_n$ for $n\geq 2$ must change its sign at least once over $\mathbb{R}$. Theorem~\ref{L:2nd} then says that each eigenfunction associated with $\lambda_2$ changes its sign either once or twice. In particular, the result is weaker than Theorem~\ref{T:sturm} when $\alpha=2$. Perhaps this is because the proof of Theorem~\ref{L:2nd} is based on a general topological argument, without taking into account the analytical properties of the operator or the eigenfunctions.  (By the way, an even eigenfunction associated with $\lambda_2$ must change its sign exactly once in the half-line $[0,\infty)$.)
\end{remark}

In the course of the proof of Theorem~\ref{L:2nd}, one establishes that if some continuous and bounded function over $\overline{\mathbb{P}^2}$ has at least two nodal domains then its trace on the boundary $\partial\mathbb{P}^2=\mathbb{R}\times \{0\}$ changes its sign at most twice. When $V$ is a periodic potential, two of the authors \cite{HJ15} extended this result and successfully developed an analogous oscillation theory for the first three periodic eigenfunctions of $H=(-d^2/dx^2)^{\alpha/2}+V(x)$, and Claassen and the second author \cite{CJ16} obtained a similar result for the first and second anti-periodic eigenfunctions. Unfortunately, it is not clear how to extend the topological arguments of \cite{FL13} and \cite{HJ15, CJ16} to higher eigenfunctions, without a messy case-by-case analysis and a strong induction argument. 

On the other hand, the key idea of the proof of Theorem~\ref{L:2nd} is to equip the totally ordered set $\{x_1^+<x_1^-<x_2^+<x_2^-\}$ with an equivalence relation $\sim$, where $x\sim y$ if $x$ and $y$ belong to the boundary of the same nodal domain of the extension of the eigenfunction to $\mathbb{P}^2$, and to bound the possible number of the equivalence classes. An appropriate combinatorial tool to substitute the above topological argument is the theory of \emph{noncrossing partitions}.

\section{Combinatorics of noncrossing partitions}\label{s:comb}

Here we develop the combinatorics of noncrossing partitions. We assume no prior familiarity with the subject on the reader's part. Throughout, the proof of Theorem~\ref{L:2nd} serves as a motivating example.

\begin{definition}
Let $X$ denote a finite and totally ordered set. A \emph{partition} $P$ of $X$ is an unordered collection of nonempty and pairwise disjoint subsets of $X$ whose union is $X$. The elements of $P$ are called \emph{blocks}. We write $|P|$ for the number of blocks in the partition $P$.
\end{definition}

A partition naturally defines an equivalence relation. Specifically, for $x,y\in X$, we write $x\sim_P y$, or $x\sim y$ for short, if $x$ and $y$ belong to the same block of the partition $P$. For $X'\subset X$, the \emph{restriction} of $P$ to $X'$ is the partition $\{B\bigcap X':B\in P\text{ and }B\bigcap X'\neq\emptyset\}$. This corresponds to the restriction of the equivalence relation $\sim_P$ to $X'$.

\begin{figure}[ht]
\begin{center}
\begin{tikzpicture}
\newcommand{\ra}{1.25} 
\newcommand{\dra}{0.1} 
\newcommand{\lra}{\ra+.4} 
\draw[thick] (0,0) circle (\ra);
\draw[ultra thick, Red, pattern = dots, pattern color = Red] (0:\ra) -- (180:\ra) -- (120:\ra) -- cycle;
\draw[ultra thick, Red] (60:\ra) -- (300:\ra);
\foreach \a in {0, 60, ..., 300}
  \draw[very thick, fill=Red] (\a:\ra) circle (\dra);
\foreach \x in {1, 2, ..., 6}
  \node at (180-\x*60:\lra) {$x_\x$}; 
\node at (-2.5,-1.5) {(a)};
\begin{scope}[shift={(7,0)}]
\draw[thick] (0,0) circle (\ra);
\draw[ultra thick, Red, pattern = dots, pattern color = Red] (0:\ra) -- (60:\ra) -- (120:\ra) -- (180:\ra) -- (300:\ra) -- cycle;
\foreach \a in {0, 60, ..., 300}
  \draw[very thick, fill=Red] (\a:\ra) circle (\dra);
\foreach \x in {1, 2, ..., 6}
  \node at (180-\x*60:\lra) {$x_\x$}; 
\node at (-2.5,-1.5) {(b)};
\end{scope}
\end{tikzpicture}
\end{center}
\caption{Configuration (a) cannot happen when the equivalence relation is ``belonging to the same nodal domain".\label{F:consolidate}}
\end{figure}
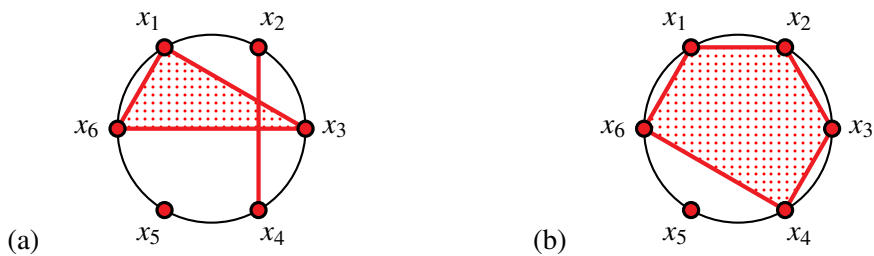

It is useful to represent a partition geometrically by placing the elements of $X$ in order around a circle, and depicting each block by a polygon whose vertices are the elements of the block. (Note that a polygon is allowed to have one or two vertices.) 
In our intended application, where the equivalence relation is ``belonging to the (boundary of the) same nodal domain",
polygons cannot cross, since nodal domains are disjoint.  For example, in Figure~\ref{F:consolidate}, if the point sets $\{x_1,x_3,x_6\}$ and $\{x_2,x_4\}$ are subsets of nodal domains $D$ and $D'$ respectively, then in fact $D=D'$, so the true equivalence relation is given by the noncrossing partition on the right.
Consequently, it suffices to work with partitions in which all associated polygons are disjoint. We make 
this precise in the following definition.

\begin{definition}
We say that two blocks $B$ and $B'$ of a partition \emph{cross} if there exist elements $x_i<x_j<x_k<x_\ell$ such that $x_i, x_k\in B$ and $x_j,x_\ell\in B'$. We say that a partition is \emph{noncrossing} if none of its blocks cross.
\end{definition}

\begin{figure}[ht]
\begin{center}
\begin{tikzpicture}
\newcommand{\ra}{1.75} 
\newcommand{\dra}{0.07} 
\newcommand{\lra}{\ra+.4} 
\draw[thick] (0,0) circle (\ra);
\draw[ultra thick, Red, pattern = dots, pattern color = Red, line join = bevel] (90:\ra) -- (157.5:\ra) -- (45:\ra) -- cycle;
\draw[ultra thick, Red, pattern = dots, pattern color = Red, line join = bevel] (202.5:\ra) -- (247.5:\ra) -- (315:\ra) -- (337.5:\ra) -- cycle;
\draw[ultra thick, Red] (22.5:\ra) -- (180:\ra);
\foreach \a in {0, 22.5, ..., 337.5}
  \draw[very thick, fill=Black] (\a:\ra) circle (\dra);
\foreach \x in {1, 2, ..., 16}
  \node at (90+360/16-\x*360/16:\lra) {$x_{\x}$}; 
\end{tikzpicture}
\end{center}
\caption{A noncrossing partition.\label{F:ncp}} 
\end{figure}
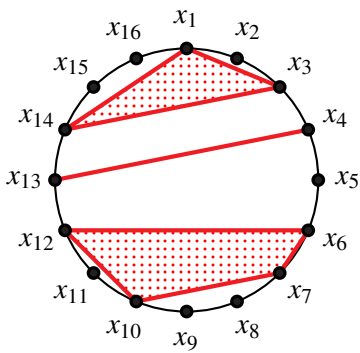

Figure~\ref{F:ncp} illustrates a noncrossing partition with blocks $\{x_1, x_3, x_{14}\}$, $\{x_4, x_{13}\}$, $\{x_6, x_7, x_{10}, x_{12}\}$, and $\{x_5\}$, $\{x_8\}$, $\{x_9\}$.
 
\begin{remark} \label{R:cyclic}
Note that the noncrossing condition is invariant under cyclic permutations of the underlying set. This will be useful later.
\end{remark}

Recall that in the course of the proof of Theorem~\ref{L:2nd}, each element of the set $\{x_1^+<x_1^-<x_2^+<x_2^-\}$ is on the boundary of either the ``positive" or ``negative" nodal domain of $E(\phi_2)$. Each nodal domain is path connected, and the positive and negative nodal domains are disjoint.  Accordingly, we may assign each element with one of two colors, depending on the sign of the nodal domain containing it, and each block must contain only points of a single color.  This motivates the following definition.

\begin{definition}[MNP]\label{def:MNP}
For a positive integer $n$, let $X^+=\{x_1^+,\dots,x_n^+\}$, $X^-=\{x_1^-,\dots,x_n^-\}$, and $X=X^+\bigcup X^-$, where $X$ is totally ordered by
\[
x_1^+<x_1^-<\cdots<x_n^+<x_n^-.
\] 
A partition $P$ of $X$ is \emph{monochromatic} of order $n$ if every block is a subset of either $X^+$ or $X^-$. We are interested in \emph{monochromatic noncrossing partitions} of $X$, or MNPs for short.  The set of all MNPs of order $n$ is denoted by $\text{MNP}_n$.
\end{definition}

\begin{figure}[ht]
\begin{center}
\begin{tikzpicture}
\newcommand{\ra}{2} 
\newcommand{\dra}{0.1} 
\newcommand{\lra}{\ra+.4} 
\draw[thick] (0,0) circle (\ra);
\draw[ultra thick, Red, pattern=grid, pattern color=Red] (0:\ra) -- (60:\ra) -- (120:\ra) -- cycle;
\draw[ultra thick, Red] (180:\ra) -- (240:\ra);
\draw[ultra thick, Blue, pattern=north east lines, pattern color=Blue] (150:\ra) -- (270:\ra) -- (330:\ra) -- cycle;
\foreach \a in {0, 60, ..., 300}
  \draw[very thick, fill=Red] (\a:\ra) circle (\dra);
\foreach \a in {30, 90, ..., 330}
  \draw[very thick, fill=Cyan] (\a:\ra) circle (\dra);
\foreach \i in {1, 2, ..., 6}
{
  \node at (150-\i*60:\lra) {$x_\i^+$}; 
  \node at (120-\i*60:\lra) {$x_\i^-$}; 
}
\end{tikzpicture}
\end{center}
\caption{A MNP of order 6.\label{F:MNP}}
\end{figure}
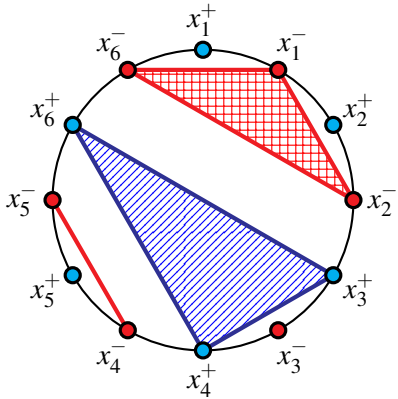

For example, $\{\{x_1^+\}, \{x_2^+\}, \{x_3^+x_4^+x_6^+\}, \{x_5^+\}, \{x_1^-x_2^-x_6^-\}, \{x_3^-\}, \{x_4^-x_5^-\}\}$ is an MNP of order $6$; see Figure~\ref{F:MNP}.

We conclude our discussion by studying the number of blocks an MNP of order $n$ may have. In the intended application, for a bounded domain $\Omega\subset\mathbb{R}^2$ with boundary $\partial\Omega$ and a continuous function on $\overline{\Omega}$, we wish to estimate the maximum number of sign changes of the function over $\partial\Omega$ in terms of an upper bound on the number of nodal domains in $\Omega$. This is equivalent to determine the minimum number of nodal domains of the function in $\Omega$ in terms of the number of sign changes on the boundary. 

\begin{lemma}[Number of blocks in an MNP] \label{P:block-size-bound}
Every MNP of order $n$ has at least $n+1$ blocks.
\end{lemma}

\begin{proof}
We use the notation of Definition~\ref{def:MNP}. Let $f(n)=\min\{|P|: P\in\text{MNP}_n\}$. Note that $f(n)\leq n+1$ for all $n$. Indeed, the partition $\{X^+,\{x_1^-\},\dots,\{x_n^-\}\}$, among others, has $n+1$ blocks.  It remains to show that $f(n)\geq n+1$.

We proceed by induction. For $n=1$, there is only one MNP of order~1, with two singleton blocks.

Let $n\geq 2$ and $P\in\text{MNP}_n$. If every block of $P$ has one element then there are $2n$ blocks and we are done. Otherwise, let $B$ denote a block with at least two elements. Without loss of generality, we may assume $B\subset X^-$.  Let $x_i^-=\min(B)$ and $x_j^-=\max(B)$, and let $k=j-i$.  Let $P'$ and $P''$ denote the restrictions of $P$ to nonempty subsets
\[
X'=\{x_1^+, x_1^-,\dots, x_i^+, x_i^-, x_{j+1}^+, x_{j+1}^-, \dots, x_n^+, x_n^-\} \quad\text{and}\quad
X''=\{x_{i+1}^+, x_{i+1}^-,\dots,x_j^+,x_j^-\},
\]
respectively. Note that $B$ is the only block of $P$ which meets both $X'$ and $X''$, since any other such block would cross $B$. Thus $|P'|+|P''|=|P|+1$, counting $B$ twice and every other block of $P$ once. Since $P'$ and $P''$ are MNP's of orders $n-k$ and $k$, respectively, it follows by induction that
\[
|P|=|P'|+|P''|-1\geq(n-k+1)+(k+1)-1=n+1.
\]
This completes the proof.
\end{proof}

\begin{remark}
Noncrossing partitions admit an involution known as \emph{Kreweras complementation} (see, e.g., \cite{Kreweras}), which answers the combinatorial question of counting MNP's and, moreover, provides a simple alternative proof of Lemma~\ref{P:block-size-bound}.  As usual, we represent a noncrossing partition $P$ by placing points labeled $1, \dots, 2n$ around a circle and associating each block with a polygon. Insert a point labeled $i'$ between the points $i$ and $i+1$ (mod $2n$) for $i=1,2,\dots, 2n$. The polygons associated with the blocks of $P$ are then regarded as to dissect the circle into cells so that every point $i'$ lies in the interior of a cell. The Kreweras complement $K(P)$ is then the partition whose blocks are the cells. For instance, if $P=\{\{2\},\{8\},\{10\},\{4,6,12\},\{1,3\},\{5\},\{7\},\{9,11\}\}$ then $K(P)=\{\{1, 2\}, \{3, 12\}, \{4, 5\}, \{6, 7, 8, 11\}, \{9, 10\}\}$; see Figure~\ref{F:mnp-enp}.

\begin{figure}[ht]
\begin{center}
\begin{tikzpicture}
\newcommand{\ra}{2.25} 
\newcommand{\dra}{0.1} 
\newcommand{\lra}{\ra+.35} 
\fill[Yellow!40!White] (0:\ra) -- (180:\ra) arc (180:360:\ra);
\fill[LimeGreen!60!White] (330:\ra) -- (270:\ra) arc (270:330:\ra);
\fill[Yellow!40!White] (120:\ra) -- (180:\ra) arc (180:120:\ra);
\fill[Yellow!40!White] (120:\ra) -- (0:\ra) arc (0:120:\ra);
\fill[LimeGreen!60!White] (30:\ra) -- (90:\ra) arc (90:30:\ra);
\draw[thick] (0,0) circle (\ra);
\draw[ultra thick] (0:\ra) -- (120:\ra) -- (180:\ra) -- cycle;
\draw[ultra thick] (30:\ra) -- (90:\ra);
\draw[ultra thick] (270:\ra) -- (330:\ra);
\foreach \a in {0, 60, ..., 300}
  \draw[very thick, fill=Cyan] (\a:\ra) circle (\dra);
\foreach \a in {30, 90, ..., 330}
  \draw[very thick, fill=Red] (\a:\ra) circle (\dra);
\foreach \a in {15, 45, ..., 345}
  \draw[thick, fill=Black] (\a:\ra) circle (\dra/2);

\foreach \i in {1, 2, ..., 12}
{
  \node at (30*\i:\lra) {\scriptsize\sf\i}; 
  \node at (15+30*\i:\lra) {\scriptsize\sf{\i}'}; 
}
\end{tikzpicture}`
\end{center}
\caption{Kreweras complementation.\label{F:mnp-enp}}
\end{figure}
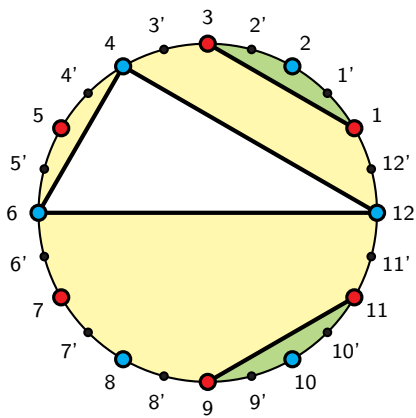

If $P$ is a noncrossing partition of $X$ with $n$ blocks, then $K(P)$ has $|X|-n+1$ blocks. Moreover, $P$ is monochromatic if and only if $K(P)$ is \emph{$2$-divisible}, i.e., all blocks have even cardinality.  This is a special case of \cite[Lemma~4.3.7]{Armstrong}, implying in particular that the number of MNP's of order $n$ is $\frac{1}{2n+1}\binom{3n}{n}$ (sequence \href{http://oeis.org/A001764}{\#A001764} in \cite{OEIS}), an instance of the \emph{Fuss-Catalan numbers} or \emph{Raney numbers} $\frac{1}{(m-1)n+1}\binom{mn}{n}$. 
Moreover, a 2-divisible partition of $\{1,\dots,2n\}$ can have at most $n$ blocks. Therefore, its Kreweras dual MNP must have at least $n+1$ blocks, offering an alternative proof of Lemma~\ref{P:block-size-bound}.
\end{remark}

\section{Applications}\label{s:application}

\subsection{Fractional Schr\"odinger operators on the real line}

We combine the results of Section~\ref{s:analysis} and Lemma~\ref{P:block-size-bound} to bound the number of sign changes in the eigenfunctions for a fractional Schr\"odinger operator on $\mathbb{R}$. 
This extends the result of \cite{FL13}
for the second eigenfunction to \emph{all} higher eigenfunctions.

\begin{theorem}[Oscillation of eigenfunctions]\label{T:line}
Let $0<\alpha<2$, $V:\mathbb{R} \to\mathbb{R}$ and $V\in L^\infty(\mathbb{R})$. For an integer $N\geq 1$,
suppose that the $L^2(\mathbb{R})$ spectrum of 
\[
H=\left(-\frac{d^2}{dx^2}\right)^{\alpha/2}+V(x)
\] 
contains at least $N$ eigenvalues
\[
\lambda_1<\lambda_2\leq\cdots\leq\lambda_N<\min\sigma_{\rm ess}(H).
\]
For each $n=1,2,\ldots, N$, if $\phi_n\in H^{\alpha/2}(\mathbb{R})\bigcap C^0(\mathbb{R})$ 
is a real eigenfunction of $H$ associated with 
the eigenvalue $\lambda_n$, then $\phi_n$ changes its sign at most $2(n-1)$ times in $\mathbb{R}$.
\end{theorem}

\begin{proof}
Suppose on the contrary that $\phi_n$ changes its sign at least $2n-1$ times in $\mathbb{R}$. Then there would exist $2n$ points in $\mathbb{R}$, say, $x_1^+<x_1^-<x_2^+<x_2^-<\cdots<x_n^+<x_n^-$ such that, without loss of generality, 
\[
\phi_n(x_k^+)>0\quad\text{and}\quad\phi_n(x_k^-)<0\qquad\text{for $k=1,2,\dots,n$}.
\] 
Note from Theorem~\ref{L:courant} that the extension $E(\phi_n)$ to the upper half plane (see \eqref{def:E}) has at most $n$ nodal domains. 

We equip the ordered set
\[
X=\{x_1^+<x_1^-<\cdots<x_n^+<x_n^-\}
\]
with an equivalence relation $\sim$, where $x\sim y$ if $x$ and $y$ belong to the boundary of the same nodal domain of $E(\phi_n)$. Let $P$ be the corresponding partition of $X$; see the previous section.  The partition $P$ is noncrossing because the nodal domains of $E(\phi_n)$ are pairwise disjoint, and it is monochromatic because the sign of $E(\phi_n)$ is constant on each nodal domain.  Thus $P\in\text{MNP}_n$, and, by Lemma~\ref{P:block-size-bound}, it has at least $n+1$ blocks, implying that $E(\phi_n)$ has at least $n+1$ nodal domains. A contradiction completes the proof.
\end{proof}

Note that the result seems \emph{not} to be sharp. Indeed, when $\alpha=2$, to compare to Theorem~\ref{T:sturm}, it appears that our method allows for a ``double counting" of the sign changes in the eigenfunctions of the classical Schr\"odinger operator. Perhaps, some analytical features of the operator $H=(-d^2/dx^2)^{\alpha/2}+V(x)$ and its eigenfunctions account for the double counting, refining the result of Theorem~\ref{T:line} to a tight bound.  
We do not pursue this here but consider it as an important question for future work.  

Note that in \cite{FL13} the result of Theorem \ref{T:line} is said to be \emph{sharp}.  There, however,
the authors seem to mean that if the eigenfunction associated with $\lambda_2$ happened to be even, then
it would indeed have exactly two sign changes over $\mathbb{R}$.  In contrast, we consider the result not to
be sharp by comparing to the classical case in Theorem~\ref{T:sturm}.

\begin{remark}\label{R:mult2}
Multiplicities in Theorem \ref{T:line} are handled in the same way as described in Remark~\ref{R:mult}.
For instance, if $\lambda_2$ has algebraic multiplicity two, so that
\[
\lambda_1<\lambda_2=\lambda_3<\lambda_4\leq\ldots\leq\lambda_N,
\]
then an eigenfunction of $H=(-d^2/dx^2)^{\alpha/2}+V(x)$ associated with $\lambda_2=\lambda_3$ has
at most two sign changes over $\mathbb{R}$, while an eigenfunction associated with $\lambda_4$ has
at most six sign changes over $\mathbb{R}$.
\end{remark}

\

The result of Theorem~\ref{T:line} may be directly adapted to radial eigenfunctions of $(-\Delta)^{\alpha/2}+V(x)$ in $\mathbb{R}^d$ for $d\geq 2$, where for $0<\alpha<2$, the fractional Laplacian $(-\Delta)^{\alpha/2}$ is defined via its multiplier $|\xi|^\alpha$ in $\mathbb{R}^d$, and $V:\mathbb{R}^d\to \mathbb{R}$ is real valued, bounded, and radial. This will extend a recent result in \cite{FLS16}. Unfortunately, the result would certainly not be sharp. Specifically, a radial eigenfunction associated with the $n$-th eigenvalue of $(-\Delta)^{\alpha/2}+V(x)$ changes its signs $2(n-1)$ times over $(0,\infty)$.  A homotopy argument in the spirit of \cite{FLS16}, or its appropriate substitute, may improve the result to the bound of $n-1$ sign changes over $(0,\infty)$, which is sharp when $\alpha=2$.  This is an interesting and important direction of future research.

\subsection{Periodic potentials}

We turn the attention to $V:\mathbb{R}\to\mathbb{R}$ smooth and $T$-periodic for some $T>0$. It is well known that the $L^2(\mathbb{R})$ spectrum of $H=(-d^2/dx^2)^{\alpha/2}+V(x)$ with a periodic potential contains no eigenvalues. Rather, the spectrum is purely essential. Moreover, if $\lambda\in\mathbb{C}$ belongs to the $L^2(\mathbb{R})$ spectrum of $H$ if and only if the quasi-periodic spectral problem
\begin{equation}\label{e:quasiper}
\begin{cases}
H\phi=\lambda \phi \qquad&\text{for $x\in(0,T)$},\\
\phi(x+T)=e^{i\xi T}\phi(x)\quad&
\end{cases}
\end{equation}
admits a nontrivial solution for some $\xi$ in the range $[-\pi/T,\pi/T)$; see, e.g., \cite{BHJ16} and references therein.  Of particular interest is when $\xi=0$ and $\xi=\pm\pi/T$, for which \eqref{e:quasiper} represents the periodic and anti-periodic eigenvalue problem, respectively.

When $\xi=0$, this amounts to studying the operator $H=(-d^2/dx^2)^{\alpha/2}+V(x)$ acting on 
\[
L^2_{\rm per}([0,T])=\{f\in L^2_{\rm loc}([0,T]): f(x+T)=f(x)\quad\text{for all $x\in\mathbb{R}$}\}
\]
with domain $H^{\alpha}_{\rm per}([0,T])$, defined analogously. Note that the $L^2_{\rm per}([0,T])$ spectrum of $H$ consists of an increasing sequence of discrete eigenvalues with finite multiplicities, accumulating only at $+\infty$. Moreover, any eigenfunction is real valued, continuous and bounded.  

\begin{theorem}[Oscillations with periodic potentials]\label{T:per}
Let $0<\alpha<2$ and $V:\mathbb{R}\to\mathbb{R}$ be smooth, bounded, and $T$-periodic for some $T>0$. Let 
\[
\lambda_1<\lambda_2\leq\lambda_3\leq\cdots\nearrow+\infty
\]
denote the $L^2_{\rm per}([0,T])$ eigenvalues of $H=(-d^2/dx^2)^{\alpha/2}+V(x)$. 
For each $n\in\mathbb{N}$, if $\phi_n\in H^{\alpha/2}_{\rm per}([0,T])\bigcap C^0([0,T])$ is a real eigenfunction of $H$ associated with the eigenvalue $\lambda_n$ then $\phi_n$ changes its sign at most $2(n-1)$ times in $[0,T)$.
\end{theorem}

The proof follows along the same line as that of Theorem~\ref{T:line}, but with appropriate modifications to accommodate the periodic boundary condition. For completeness, we sketch the main ideas of the proof.

First, note that the fractional Laplacian $(-d^2/dx^2)^{\alpha/2}$ acting on $L^2_{\rm per}([0,T])$ may be viewed as the Dirichlet-to-Neumann operator for an appropriate local problem in the periodic half strip $[0,T]_{\rm per}\times(0,\infty)$; see \cite{RS16} and \cite{HJ15} for details.  
We then derive a variational characterization of the eigenvalues and eigenfunctions for $H=(-d^2/dx^2)^{\alpha/2}+V(x)$ in terms of the Dirichlet type functional 
\[
\iint_{(0,T)\times(0,\infty)}y^{1-\alpha}|\nabla w(x,y)|^2~dx\;dy+\int_0^T V(x)|w(x,0)|^2~dx,
\]
which is defined in an appropriate function space for $w=w(x,y)$ in $[0,T]_{\rm per}\times(0,\infty)$, where $w(x,0)$ denotes the trace of $w(x,y)$ on the boundary $[0,T]_{\rm per}\times \{0\}$. 
An argument in the spirit of Courant's nodal domain theorem implies that the extension of $\phi_n$ to the periodic half strip $[0,T]_{\rm per}\times(0,\infty)$ has at most $n$ nodal domains. We merely pause to remark that if an eigenfunction of $H$ is bounded and continuous in $[0,T]$ then its extension to the periodic half strip belongs to $C^0([0,T)_{\rm per}\times(0,\infty))$. Finally, we combine this and Lemma~\ref{P:block-size-bound} to complete the proof. We omit the details.

Theorem~\ref{T:per} extends the result in \cite{BH14, HJ15} for $\phi_2$ and $\phi_3$ to all higher eigenfunctions, which by the way plays a central role in the study of the stability and instability of periodic traveling waves of fractional Korteweg-de Vries equations, with respect to period preserving perturbations as well as long wavelength perturbations. 

Theorem~\ref{T:per} agrees with the Sturm oscillation theorem when $\alpha=2$; see, e.g., \cite{KPbook}. 
Like Theorem \ref{T:line}, it is not sharp, as the following elementary
example shows.

\begin{example}\label{ex:per}
When $V=0$, so that $H=(-d^2/dx^2)^{\alpha/2}$ is a fractional Laplacian, a straightforward calculation reveals that the $L^2_{\rm per}([0,2\pi])$ eigenvalues of $H$ are 
\[
0<1=1<2=2<3=3<\ldots\nearrow+\infty.
\]
The lowest eigenvalue $\lambda_1=0$ is simple and an associated eigenfunction is $1$, which is sign definite.  All
higher eigenvalues have algebraic multiplicity two. For each $n\in\mathbb{N}$, the eigenvalues
$\lambda_{2n}=\lambda_{2n+1}=n$ have the eigenspace spanned by $\cos(nx)$ and 
$\sin(nx)$, which change their signs exactly $2n$ times over $[0,2\pi)$.  In particular,
recalling how multiplicities are handled (see Remark \ref{R:mult2}), we find that Theorem~\ref{T:per} 
over counts the oscillation by $n-2$, where $n$ denotes the \emph{lowest} index
for which the eigenvalue is repeated.

Interestingly, the over counting in the example seems to stem from that non-zero eigenvalues
have algebraic multiplicity two. In fact, if we restrict to the sector of $L^2_{\rm per}([0,2\pi])$
of even functions, then all eigenvalues are simple and the bound in Theorem~\ref{T:per} is sharp.
Perhaps, this should not be surprising considering that the way that the Courant's nodal domain theorem handles eigenvalues with multiplicities is not sharp.
\end{example}

The above may possibly apply to other boundary conditions, for instance, for a fractional Schr\"odinger operator with a periodic potential acting on a space of ``anti-periodic" functions, as considered in \cite{CJ16}, corresponding to 
$\xi=\pm\pi/T$ in \eqref{e:quasiper}.  This is an interesting direction of future research.

Our methods should be directly applicable to general pseudo-differential operators beyond the fractional Laplacian, provided that they may be regarded as Dirichlet-to-Neumann
operators for appropriate, local elliptic problems. 
Examples include $(-d^2/dx^2+m^2)^{\alpha/2}$, for $\alpha\in(0,2)$ and $m>0$, in the theory of magnetic Schr\"odinger operators,
and $(-d^2/dx^2)^{1/2}\coth(-d^2/dx^2)^{1/2}$ in the intermediate long-wave equation; see, e.g., \cite{FLS16}.

\subsection{Steklov problems} \label{s:steklov}

Let $\Omega\subset\mathbb{R}^2$ be a smooth bounded domain, i.e., an open and connected set whose closure $\overline{\Omega}$ 
is compact with smooth boundary $\partial\Omega$.
The Steklov problem on $\Omega$ is 
\begin{equation}\label{e:steklov}
\begin{cases}
\Delta u=0 \qquad & \textrm{in}\quad \Omega,\\
\frac{\partial u}{\partial n}=\lambda u & \textrm{on}\quad \partial\Omega,
\end{cases}
\end{equation}
where $\Delta$ denotes the Laplacian in $\mathbb{R}^2$ and $\partial/\partial n$ is the outward normal derivative along the boundary $\partial\Omega$. 

This problem was introduced by the Russian mathematician V. A. Steklov at the turn of the last century. We refer the interested reader to \cite{GP14} for a recent survey. It is well known (see, e.g., \cite{KS69}) that the spectrum of \eqref{e:steklov} consists of an increasing sequence of discrete eigenvalues, accumulating only at $+\infty$, denoted
\[
0=\lambda_1<\lambda_2\leq\lambda_3\leq\cdots\nearrow+\infty.
\]
It is straightforward to verify that the lowest eigenvalue $\lambda_1$ is simple.
Higher eigenvalues may have nontrivial multiplicities.

The Steklov eigenvalues may be interpreted as the eigenvalues 
of the Dirichlet-to-Neumann operator which maps a function $f$ on $\partial\Omega$ to $\partial E(f)/\partial n$, 
where $E(f)$ is the harmonic extension of $f$ to $\Omega$. 
Consequently, it seems natural to expect that the methodologies presented here may be applicable
to \eqref{e:steklov}.  The goal of this section is to carry out the program.

To begin, we observe that Courant's nodal domain theorem applies to the Steklov problems mutatis mutandis.
See, e.g., \cite{KKP14} for details. 

\begin{lemma}\label{l:steklovnd}
Let $\Omega$ be a smooth bounded domain in $\mathbb{R}^2$. Then a Steklov eigenfunction associated with the eigenvalue $\lambda_n$ has at most $n$ nodal domains in $\Omega$.
\end{lemma}

The next result is an analogue of Lemma~\ref{P:block-size-bound} for bounded domains $\Omega\subset\mathbb{R}^2$ of arbitrary genera.  The argument is topological rather than analytical: $\Omega$ need not be smooth, and the function $f:\overline{\Omega}\to\mathbb{R}$ whose sign changes and nodal domains we study need only be continuous.

\begin{theorem}
\label{bound-general}
Let $\Omega$ be a bounded domain in $\mathbb{R}^2$ of genus $g$ such that the components of $\partial\Omega$ consist of simple closed curves $\gamma_0,\dots,\gamma_g$.
Let $f:\overline{\Omega}\to\mathbb{R}$ be a continuous function. Then
\begin{equation} \label{most-general-bound}
s\leq 2(n+g-1),
\end{equation}
where $n$ is the number of nodal domains of $f$ in $\Omega$ and $s$ is the number of sign changes of $f$ along $\partial\Omega$.  
\end{theorem}

\begin{proof}
We proceed by induction on $g$. First, suppose $g=0$; note that $\gamma_0=\partial\Omega$. Suppose that $f$ changes its sign at least $2n-1$ times on $\partial\Omega$. Then there exist points $x_1^+,x_1^-,x_2^+,x_2^-,\dots,x_n^+,x_n^-$, listed in cyclic order along $\partial\Omega$, such that $f(x_k^+)>0$ and $f(x_k^-)<0$ for $k=1,2,\dots, n$.
The argument of Lemma~\ref{P:block-size-bound} applies verbatim to prove the assertion. Note that the choice of $x_1^+$ does not matter, since the monochromatic property is invariant under cyclic permutations of the underlying set; see Remark~\ref{R:cyclic}. 

Now suppose $g>0$. For each $0\leq i\leq g$, let $s_i$ be the number of sign changes on the curve $\gamma_i$.
By the genus-zero argument, each $\gamma_i$ gives rise to at least $s_i/2+1$ nodal domains.  If no nodal domain meets more than one $\gamma_i$, then $n\geq\sum_{i=0}^g(s_i/2+1) = s/2+g+1$ or, equivalently, $s\leq 2(n-g-1)$.
Otherwise, suppose that $x$ and $x'$ are points on different boundary components $\gamma$ and $\gamma'$ such that $x$ and $x'$ belong to the same nodal domain $D$, which may be assumed to be positive.  Let $C$ be a continuous simple path in $D$ connecting $x$ to $x'$; see Figure~\ref{F:high-genus}.

\begin{figure}[h]
\begin{center}
\resizebox{4in}{2in}{\includegraphics{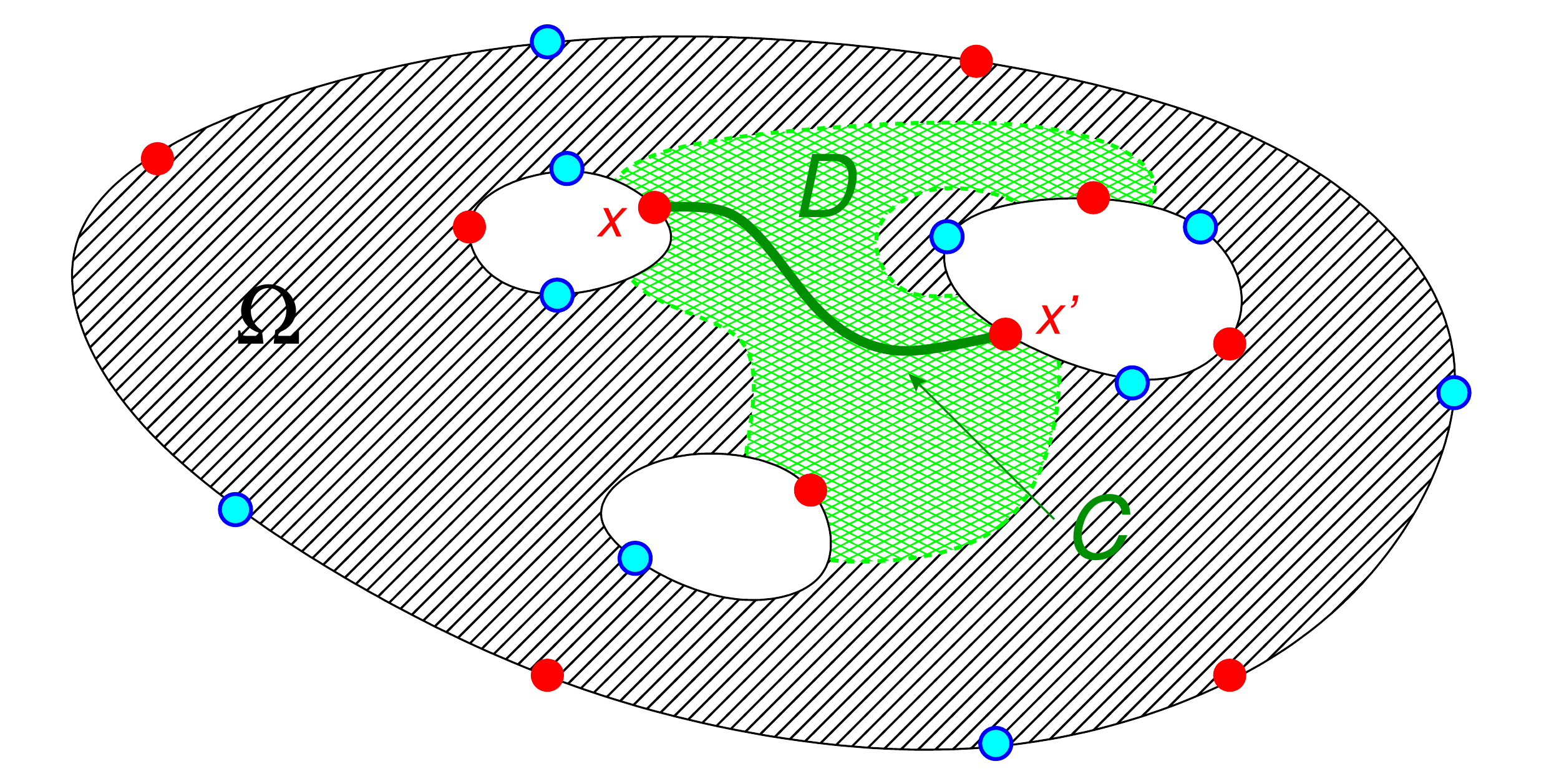}}
\end{center}
\caption{A domain $\Omega\subset\mathbb{R}^2$ of high genus with sign changes along boundary curves. \label{F:high-genus}}
\end{figure}

First, we show that we may assume without loss of generality that $f$ is constant on some tubular neighborhood $U$ of $C$.
Indeed, let $U$ and $V$ be tubular neighborhoods of $C$ such that $\overline{U}\subset V$ and $\overline{V}\subset D$.  By linear interpolation on the fibers of $V-U$, we can define a continuous function $\widetilde{f}:\Omega\to\mathbb{R}$ which is a positive constant $K$ on $U$ and coincides with $f$ on $\Omega-V$.  In particular, $\widetilde{f}$ has the same nodal domain structure as $f$, and we may replace $f$ by $\widetilde{f}$.

Second, we deform the pair $(\Omega,f)$ continuously by shrinking the curve $C$ to a point.  To make this construction precise, note that we may assume in addition that the tubular neighborhood $U\supset C$ is thin enough that it avoids all components of $\partial\Omega$ other than $\gamma$ and $\gamma'$.  Thus we have a homeomorphism $\varphi:U\to U'_0=[-1,1]\times[-1,1]$ such that $\varphi(C)=[-1,1]\times\{0\}$, and $\{-1\}\times[-1,1]$ and $\{1\}\times[-1,1]$ correspond to sections of $\gamma$ and $\gamma'$, respectively; see Figure~\ref{F:tubular}.

\begin{figure}[h] 
\begin{center}
\resizebox{6in}{1.3125in}{\includegraphics{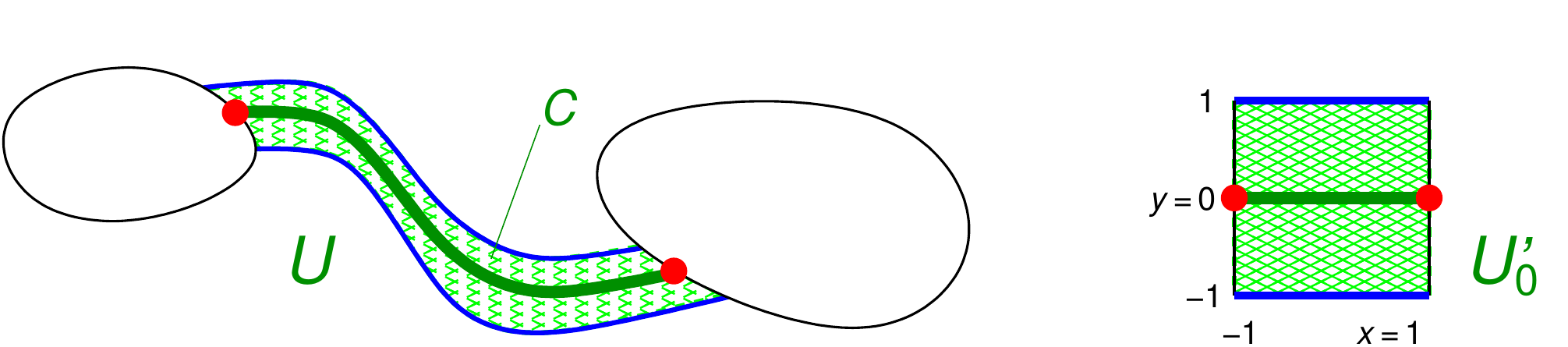}} 
\end{center}
\caption{Coordinatizing the tubular neighborhood $U$.\label{F:tubular}}
\end{figure}

\noindent
For $0\leq t\leq 1$, let
\[U'_t=\{(x,y):\ |y|\leq 1,\ |x|\leq 1-t(1-|y|)\}\quad\text{and}\quad U_t=\varphi^{-1}(U'_t)\]
(see Figure~\ref{F:shrink}) and let $\Omega_t=(\Omega\setminus U)\bigcup U_t$.  In addition, define $f_t:\Omega_t\to\mathbb{R}$ by $f_t|_{\Omega\setminus U}=f|_{\Omega\setminus U}$ and $f_t|_{U_t}\equiv K$, so that $f_t$ is continuous on $\Omega_t$ for all $t$.  The domains $\Omega_t$ are mutually homeomorphic for $0\leq t<1$, and the numbers of sign changes on each boundary curve and the number of nodal domains remain constant in the range.

\begin{figure}[h] 
\begin{center}
\resizebox{6in}{.875in}{\includegraphics{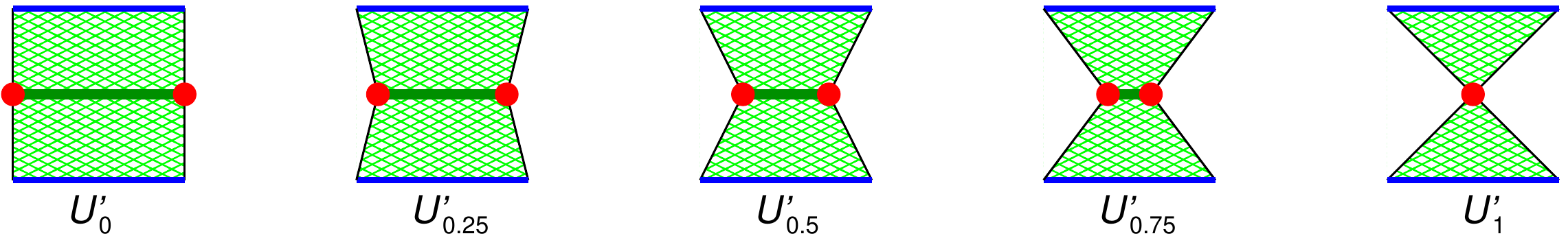}} 
\end{center}
\caption{The sets $U'_t$.  Note that the top and bottom parts of the boundary are constant for all $t$.\label{F:shrink}}
\end{figure}

When $t=1$, we obtain a space $\Omega_1$ in which the two simple curves $\gamma$ and $\gamma'$ have merged into a single curve with a node $z$.  This curve can be perturbed into a simple closed curve, re-splitting the double point into two points $z'$ and $z''$ to obtain a domain $\Omega'$ as in Figure~\ref{F:deform}.  This splitting gives rise to a continuous function $\rho:\Omega'\to\Omega_1$ that maps each of $z'$ and $z''$ to $z$ and is a homeomorphism elsewhere, and we can define $f':\Omega'\to\mathbb{R}$ by pullback as $f'=f_1\circ\rho$.

\begin{figure}[h]
\begin{center}
\resizebox{6.075in}{1.425in}{\includegraphics{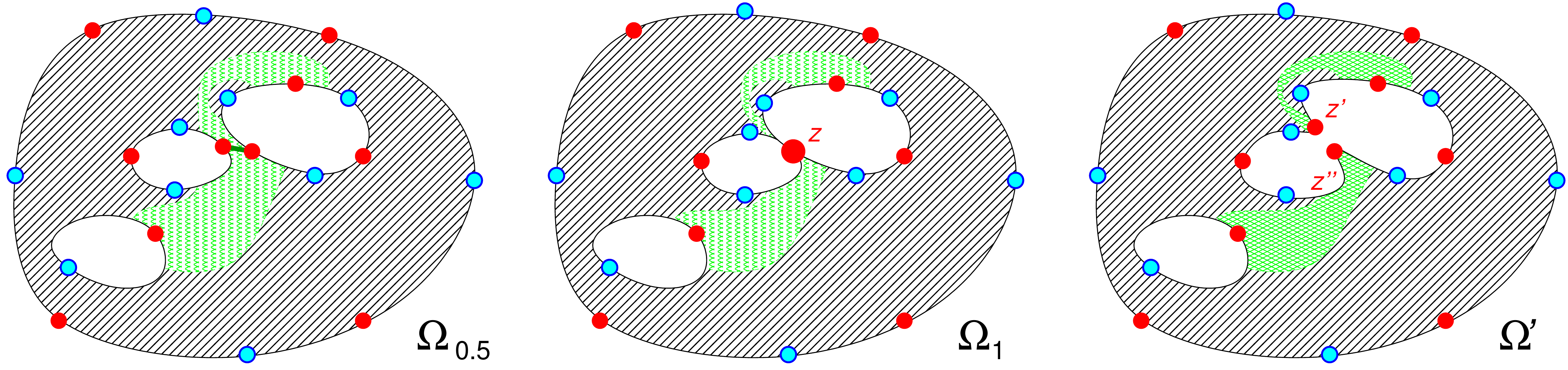}} 
\end{center}
\caption{Deformation of $\Omega$ into a domain of lower genus. \label{F:deform}}
\end{figure}

By construction, the domain $\Omega'$ has genus $g-1$.  Moreover, the continuous function $f':\Omega'\to\mathbb{R}$ has the same number of sign changes as $f$, and either $n$ or $n+1$ nodal domains (the latter if and only if the original nodal domain $D$ was simply connected, as in the example).  By induction on $g$, we obtain 
\[
s\leq 2((n+1)+(g-1)-1) =  2(n+g-1).
\]
This completes the proof.
\end{proof}

It is worth emphasizing that the result holds for all continuous functions on compact domains in $\mathbb{R}^2$; no stronger hypotheses on $f$ are required.  Moreover, the proof shows that the maximum number of sign changes is achieved when every nodal domain is simply connected.  It should be possible to obtain a tighter bound in terms of the genera of the individual nodal domains, using similar topological and combinatorial arguments, but we do not pursue this idea here.

We now apply these topological facts to the Steklov problem.

\begin{theorem}\label{t:steklov}
Under the hypothesis of Lemma~\ref{l:steklovnd}, a Steklov eigenfunction associated with the Steklov eigenvalue $\lambda_n$ changes its sign at most $2(n+g-1)$ times on $\partial\Omega$, where $g$ is the genus of $\Omega$.
\end{theorem}

We omit the details.  While the result of Theorem \ref{t:steklov} is not sharp in the case of genus zero,
interestingly, we find that it \emph{is sharp} in the case of genus $1$.  We illustrate this
in the following examples.

\begin{example}
Let $\Omega$ be the unit disk in $\mathbb{R}^2$. Note that $\Omega$ is simply connected. It is straightforward that the Steklov eigenvalues are 
\[
0<1=1<2=2<3=3<\dots\nearrow+\infty.
\]
The lowest eigenvalue $\lambda_1=0$ is simple with eigenfunction $1$, which is sign definite.
All other higher eigenvalues have algebraic multiplicity two. For each $n\in\mathbb{N}$, the eigenvalues
$\lambda_{2n}=\lambda_{2n+1}=n$ have the eigenspace spanned by $r^n\cos(n\theta)$ and $r^n\sin(n\theta)$ in polar coordinates.  Restricting to the boundary $r=1$, the eigenfunctions agree with those in Example~\ref{ex:per},
implying that the result of Theorem~\ref{t:steklov} is not sharp.
\end{example}

\begin{example}
For $\varepsilon\in(0,1)$, let $\Omega=B(0,1)\setminus \overline{B(0,\varepsilon)}$ be the annulus consisting of the open unit disk in $\mathbb{R}^2$ with the closure of the disk $B(0,\varepsilon )$ 
removed. It is known that the lowest eigenvalue
\[
\lambda_1=\frac{1+\varepsilon}{\varepsilon\ln(1/\varepsilon)}
\] 
is simple with a radially symmetric and sign definite eigenfunction
\[
\phi_1(r)=-\left(\frac{1+\varepsilon}{\varepsilon\ln(\varepsilon)}\right)\ln(r)+1
\]
in polar coordinates.
Higher eigenvalues may be ordered as
\[
\lambda_1<\lambda_2=\lambda_3<\lambda_4=\lambda_5<\dots\nearrow+\infty;
\]
see \cite{GP14} for details.
For each $n\in\mathbb{N}$, the eigenvalues $\lambda_{2n}=\lambda_{2n+1}$ have the eigenspace spanned by
\[
(A_nr^n+B_n r^{-n})\cos(n\theta)\quad{\rm and}\quad (A_nr^n+B_n r^{-n})\sin(n\theta)
\]
in polar coordinates
for appropriate real constants $A_n$ and $B_n$.  It is then easy to verify that
the eigenfunctions associated to $\lambda_{2n}=\lambda_{2n+1}$ changes their signs exactly $4n$ times on the boundary. Recalling
how multiplicities are handled (see Remark~\ref{R:mult2}), it follows that
the oscillation
bound in Theorem~\ref{t:steklov} is achieved in this example, whence it is sharp. 
\end{example}

\section*{Acknowledgments} 
The authors thank Jared Bronski, Graham Cox, Rick Laugessen, Saul Stahl,
Selim Sukhtaiev, and Alim Sukhtayev for helpful conversations, and the anonymous referee for several useful comments. 
VMH thanks the Mathematics Department at Brown University for its generous hospitality. 

\bibliographystyle{amsplain}


\begin{dajauthors}
\begin{authorinfo}[vmh]
  Vera Mikyoung Hur\\
  University of Illinois\\
  Urbana, Illinois, USA\\
  verahur\imageat{}math\imagedot{}uiuc\imagedot{}edu\\
  \url{http://www.math.illinois.edu/~verahur/}
\end{authorinfo}
\begin{authorinfo}[maj]
  Mathew A. Johnson\\
  University of Kansas\\
  Lawrence, Kansas, USA\\
  matjohn\imageat{}ku\imagedot{}edu\\
  \url{http://www.people.ku.edu/~m079j743/}
\end{authorinfo}
\begin{authorinfo}[jlm]
  Jeremy L. Martin\\
  University of Kansas\\
  Lawrence, Kansas, USA\\
  jlmartin\imageat{}ku\imagedot{}edu\\
  \url{http://www.people.ku.edu/~jlmartin/}
\end{authorinfo}
\end{dajauthors}

\end{document}